\documentclass[11pt]{amsart}
\usepackage{amssymb,mathrsfs,color,mathtools}
\usepackage[all]{xy}
\usepackage{float}
\usepackage[x11names]{xcolor}
\usepackage{tkz-graph}
\usetikzlibrary{decorations.markings}
\tikzstyle{vertex}=[circle, draw, inner sep=0pt, minimum size=6pt]

\newcommand{\luk}{\L u\-ka\-si\-e\-w\-icz}
\newtheorem{theorem}{Theorem}[section]
\newtheorem{lemma}[theorem]{Lemma}

\newtheorem{proposition}[theorem]{Proposition}

\theoremstyle{definition}
\newtheorem{definition}[theorem]{Definition}

\newtheorem{example}[theorem]{Example}

\newenvironment{proof2}
{\begin{trivlist} \item[] {\em Proof. }}
{
\hfill$\dashv$\end{trivlist}}

\newcommand{\upminus}{\rule[2.7pt]{5pt}{.4pt}}
\newcommand{\mindiamond}{\mathord{\lozenge\hspace{-5.8pt}\footnotesize{\upminus}}\,}
\newcommand{\minbox}{\boxminus}

\newcommand{\comment}[1]{}

\begin{document}

\title{Rotations of G\"odel algebras with modal operators 
}
%
%
\author[]{Tommaso Flaminio$^1$ \and
Lluis Godo$^1$\and Paula Mench\'on$^2$ \and Ricardo O. Rodriguez$^3$\\
$^1${\tiny Artificial Intelligence Research Institute (IIIA - CSIC), Barcelona, Spain. \email{\{tommaso,godo\}@iiia.csic.es}} \\
$^2${\tiny Nicolaus Copernicus University, Toru\'n, Poland, CONICET, Univ. Nacional del Centro de la Provincia de Buenos Aires, Tandil, Argentina. \email{paulamenchon@gmail.com}} \\
$^3$ {\tiny UBA-FCEyN, Departamento de Computaci\'on,\\
 CONICET-UBA, Inst. de Invest. en Cs. de la Computaci\'on, Buenos Aires, Argentina. \email{ricardo@dc.uba.ar}}\\
}
%
%

%
\maketitle   

\begin{abstract}
The present paper is devoted to study the effect of connected and disconnected rotations of G\"odel algebras with operators grounded on directly indecomposable structures. The structures resulting from this construction we will present are nilpotent minimum (with or without negation fixpoint, depending on whether the rotation is connected or disconnected) with special modal operators defined on a directly indecomposable algebra. In this paper we will present a (quasi-)equational definition of these latter structures. Our main results show that directly indecomposable  nilpotent minimum algebras (with or without negation fixpoint) with modal operators are fully characterized as connected and disconnected rotations of directly indecomposable G\"odel algebras endowed with modal operators. 
\end{abstract}

\section{Introduction}

Fuzzy modal logic is an active and rapidly growing area of research that aims at generalizing classical modal logic to a many-valued or fuzzy framework. Although the birth of Mathematical Fuzzy logic as a discipline is usually taken with the publication of H\'ajek's book \cite{H98} in 1998, the first attempts to generalize modal logic to the setting of many-valued logic can be traced back to the 90's of the last century. Indeed, in 1991 and 1992 Fitting published two fundamental papers  \cite{Fitting1,Fitting2} in which, he investigates two families of many-valued modal logics: on the side of their relational semantics, the first one is characterized by Kripke models in which, at each possible world, formulas are evaluated by a finite Heyting algebra; the second one allows also the accessibility relation to be many-valued.

Fitting's work on many-valued modal logics paved the way to the birth of fuzzy modal logic as a field and inspired several other researchers who, following his ideas, further generalized classical modal logic to the ground of infinite-valued fuzzy logics. In this respect it is worth to remember the general approaches collected in the papers by Priest \cite{Priest}, by H\'ajek \cite{HajekS5}, by Bou, Esteva, Godo and Rodriguez  \cite{BEGR} and by Diaconescu, Metcalfe and Schn\"uriger \cite{DMS}. Further works which follow Fitting's ideas of generalizing Kripke models to the fuzzy environment,  but focus on modal expansions of specific propositional  logics, are also worth to recall. In particular the ones by Caicedo and Rodriguez \cite{XARO1,XARO} who investigate modal expansions of G\"odel logic, the paper \cite{HanTe} by Hansoul and Teheux who instead consider modal expansions of propositional \luk\ logic, and the one by Vidal, Esteva and Godo \cite{VEG17} where the propositional base is product fuzzy logic.

In the recent paper \cite{FGR19}, taking inspiration on intuitionistic modal logics, we made the first steps towards an algebraic approach to modal expansions of G\"odel propositional logic by introducing  {\em finite G\"odel algebras with operators} (GAOs for short) and a class of relational models which, in contrast to the original Fitting's approach, are not Kripke-like structures. Indeed, the relational duals of GAOs turn out to be  frames based on {\em forests} rather than sets, that is, posets in which the downset of each element is totally ordered. 



In the present paper, we are interested in studying {\em nilpotent minimum algebras with modal operators}.  Nilpotent minimum algebras (NM-algebras for short) are the algebraic semantics of the so-called nilpotent minimum fuzzy logic, denoted by NM, which turns out to be sound and complete w.r.t. the algebra on $[0, 1]$ defined by the so-called nilpotent minimum t-norm and its residuum  \cite{EG01}, that is a left-continuous, but not continuous, t-norm. 
The results on NM-algebras with modal operators are obtained by extending the well-known construction of connected and disconnected rotations of G\"odel algebras \cite{Jen,Busaniche} to the modal framework. 

This paper is organized as follows.  Section \ref{sec:intro} is devoted to recall basic algebraic notions and, in particular, to introduce G\"odel and NM-algebras. Also we  recall how NM-algebras with or without negation fixpoint can be constructed as connected or disconnected rotations of  directly indecomposable G\"odel algebras. In Section \ref{sec:NMOperators}, after briefly recalling G\"odel algebras with modal operators, we study modal expansions of NM$^+$ and NM$^-$-algebras obtained by extending connected and disconnected rotations to the modal setting. In particular,  Subsection \ref{secModalNM+} is on modal NM$^+$-algebras, while NM$^-$-algebras are the subject of Subsection \ref{secModalNM-}. Our main results, namely Theorems \ref{thm:Gao-NM+} and \ref{thm:Gao-NM-}, show that for every  G\"odel algebra with operators whose G\"odel reduct is directly indecomposable and  satisfying certain properties, one can build both an NM$^+$ and an NM$^-$-algebra with operators whose NM-reduct still is directly indecomposable. Indeed, each such modal NM-algebra arises in those ways. We conclude with Section \ref{secFinal} where we present some conclusions and  some prospects for our future work.

\section{G\"odel and nilpotent minimum algebras}\label{sec:intro}
Most of the algebraic structures we consider in this paper lay in the variety of MTL-algebras \cite{EG01}, namely, integral, commutative, bounded and prelinear residuated lattices which are structures of the form ${\bf A}=(A, \ast, \to, \wedge, \vee, \bot, \top)$ of type $(2,2,2,2,0,0)$ where: $(A, \ast, \top)$ is a commutative monoid; $(A,\wedge, \vee, \bot,\top)$ is a bounded lattice; and the following conditions hold for all $x,y,z\in A$:
\begin{description}
\item[${\bf (Res)}$] $x\ast y\leq z \mbox{ iff } x\leq y\to z$, where $\leq$ denotes the lattice order of ${\bf A}$;
\item[${\bf (Pre)}$] $(x\to y)\vee (y\to x)=\top.$
\end{description}
In every MTL-algebra,  a negation operator $\neg$ can be defined as: $\neg x=x\to \bot$.

{
Let ${\bf A}$ be an MTL-algebra. A subset $f$ of ${\bf A}$ is said to be a {\em filter} provided that: (1) $\top\in f$, (2) if $x,y\in f$, then $x\ast y\in f$, (3) if $x\in f$ and $y\geq x$ then $y\in f$. A filter $f\neq A$ (that is, a {\em proper} filter) is said to be {\em prime} if $x\vee y\in f$ implies that either $x\in f$ or $y\in f$. A filter $f$ is {\em principal} (or {\em principally generated}) if there exists an element $x\in A$ such that $f={\uparrow}x=\{y\in A\mid y\geq x\}$. 

An MTL-algebra ${\bf A}$ is said to be {\em directly indecomposable} (d.i., henceforth) if it cannot be factorized as a non-trivial direct product $\prod_{i\in I}{\bf A}_i$ of MTL-algebras ${\bf A}_i$. 
}

\begin{definition}
A {\em G\"odel algebra} is MTL-algebra
which further satisfies the following equation which expresses the {\em idempotency property}:
\begin{description}
\item[${\bf (Idem)}$] $x\ast x=x$.
\end{description}
\end{definition}
The idempotency of G\"odel algebras allows us to present them in the simplified signature in which $\ast=\wedge$, which coincides with the signature of Heyting algebras \cite{Hasimoto2001,Horn}. Indeed, an equivalent definition of G\"odel algebras is to present them as those Heyting algebras, i.e., integral, commutative, bounded and idempotent residuated lattices, that further satisfy the prelinearity equation ${\bf (Pre)}$ above. 

{A G\"odel algebra is d.i. iff  $\bot$ is meet irreducible. Equivalently a  G\"odel algebra is d.i. iff it has a join-irreducible co-atom and hence a unique maximal filter that is in fact principally generated by its unique co-atom.} 
%
%
%
In what follows we will need the following easy result.

\begin{lemma}\label{GodProp1}
In every d.i.  G\"odel algebra, if $x>\bot$, $\neg x=\bot$. 
\end{lemma}
\begin{definition}
A {\em nilpotent minimum} algebra ${\bf A}$ is 
a MTL-algebra
which further satisfies the following equations:
\begin{description}
\item[${\bf (Inv)}$] $\neg\neg x=x$,
\item[${\bf (NM)}$] $\neg(x\ast y)\vee (x\wedge y\to x\ast y)=\bot$. 
\end{description}
\end{definition}
{In linearly ordered nilpotent minimum algebras, condition (NM) implies that either $x * y = \bot$ or $x * y = \min(x, y)$}. Since the negation operator of NM-algebras is involutive, if it has a fixpoint, it is unique. The property of having a negation fixpoint or not, can be captured equationally. An NM-algebra ${\bf A}$ is said to be {\em without negation fix-point} (and we write that ${\bf A}$ is an NM$^-$-algebra) if it satisfies the further equation:
\begin{description}
\item[${\bf (NM^-)}$] $\neg(\neg (x\ast x))\ast\neg(x\ast x)=(\neg(\neg x\ast \neg x))\ast(\neg(\neg x\ast\neg x))$.
\end{description}
As for those NM-algebras with a negation fixpoint, we need to expand their language by a new constant ${\bf f}$ and an algebra ${\bf A}=(A, \ast, \to, \wedge, \vee,{\bf f}, \bot, \top)$ is said to be {\em with negation fix-point} (and we write that ${\bf A}$ is an NM$^+$-algebra) if its $\{{\bf f}\}$-free reduct is a NM-algebra and it further satisfies:
\begin{description}
\item[${\bf (NM^+)}$] $\neg{\bf f}={\bf f}$.
\end{description}
{
Similarly to G\"odel algebras, {a  NM$^+$ or NM$^-$-algebra is directly indecomposable} iff it has a unique maximal filter. In this case, however, the unique maximal filter is not principally generated by a unique co-atom of ${\bf A}$. Indeed, as we will see below, there exist d.i. NM$^+$ and NM$^-$-algebras without a unique co-atom.}

{Now we  recall a general construction which allows to define NM$^-$ and NM$^+$-algebras as disconnected and connected rotations of a directly indecomposable G\"odel algebra.
The original ideas are from \cite{Jen,Busaniche} while the proofs of next results, namely Propositions \ref{propFix} and \ref{prop:NMtoGodel} and Theorems \ref{thm:isoG-NM} and \ref{thm:DINMandG}, can be found in \cite[\S4]{Busaniche}.}

Let ${\bf A}$ be a G\"odel algebra and define
\begin{description}
\item $
NM^{+}({\bf A})=\{(a^-, a^+)\in A\times A\mid a^-\wedge a^+=\bot\}
$
\item 
$NM^{-}({\bf A})=\{(a^-, a^+)\in A\times A\mid (a^-\wedge a^+)\vee \neg(a^-\vee a^+)=\bot\}. 
$
\end{description}
Further, for every $(a^-, a^+), (b^-, b^+)\in NM^+(A) \cup NM^-(A)$, define:
\begin{description}
\item $(a^-, a^+) \ast (b^-, b^+)=((a^+\vee b^+)\to(a^-\vee b^-), a^+\wedge b^+)$;
\item $(a^-, a^+)\wedge (b^-, b^+)=(a^-\vee b^-, a^+\wedge b^+)$;
\item $\neg(a^-, a^+)=(a^+, a^-)$.
\end{description}
Now, denote by ${\bf NM}^+({\bf A})$ and ${\bf NM}^-({\bf A})$ the structures:
\begin{description}
\item $
{\bf NM}^+({\bf A})=(NM^+({\bf A}), \ast, \wedge, \neg, (\bot,\bot), (\top, \bot), (\bot,\top)),
$
\item$
{\bf NM}^-({\bf A})=(NM^-({\bf A}), \ast, \wedge, \neg, (\top, \bot),(\bot,\top)).
$
\end{description}
It is worth pointing out that, in case ${\bf A}$ is a d.i.  G\"odel algebra, the above constructions that build the algebras ${\bf NM}^+({\bf A})$ and ${\bf NM}^-({\bf A})$ coincide, respectively, with the well-known {\em connected} and {\em disconnected} rotations of ${\bf A}$ (see \cite{Jen} and \cite{Busaniche}). 

\begin{proposition}\label{propFix}
For every  G\"odel algebra ${\bf A}$,  ${\bf NM}^+({\bf A})$ is a NM$^+$-algebra with negation fixpoint $(\bot, \bot)$ and  ${\bf NM}^-({\bf A})$ is a NM$^-$-algebra without negation fixpoint. Moreover, ${\bf A}$ is d.i.  iff so are ${\bf NM}^+({\bf A})$ and ${\bf NM}^-({\bf A})$. 
\end{proposition}
According to the above construction, if $(a^-, a^+), (b^-, b^+)$ are elements of either $NM^+({\bf A})$ or $NM^-({\bf A})$, then $(a^-, a^+)\leq (b^-, b^+)$ iff $(a^-, a^+)=(a^-, a^+)\wedge (b^-, b^+)=(a^-\vee b^-, a^+\wedge b^+)$ iff $a^-\geq b^-$ and $a^+\leq b^+$. 

Let us now start from any NM$^+$ or NM$^-$-algebra $\bf B$ and define
$$
G({\bf B})=\{b\ast b\mid b\in B\}. 
$$
\begin{proposition}\label{prop:NMtoGodel}
For every  NM$^+$ or $NM^-$-algebra ${\bf B}$, the structure ${\bf G}({\bf B})=(G({\bf B}), \wedge, \to^2, \bot, \top)$, where for every $x, y\in G({\bf B})$, $x\to^2 y=(x\to y)\ast(x\to y)$, is a G\"odel algebra. Further, ${\bf B}$ is d.i.  iff so is ${\bf G}({\bf B})$.
\end{proposition}
Let ${\bf A}$ be a d.i.  G\"odel algebra and let ${\bf B}$ be a d.i.  NM-algebra. Consider the  maps $\gamma: {\bf A}\to {\bf G}({\bf NM}^{\pm}(\bf A))$ and $\eta: {\bf B}\to{\bf NM^{\pm}(\bf G({\bf B}))}$ defined by the following stipulations: for all $a\in A$ and $b\in B$,
\begin{equation}\label{gamma}
\gamma(a)=(\bot, a)\ast(\bot, a),
\end{equation}
and
\begin{equation}\label{eta}
\eta(b)=\left\{
\begin{array}{ll}
(\bot, b\ast b)&\mbox{ if }b>\neg b;\\
(\neg b\ast \neg b, \bot)&\mbox{ if }b\leq \neg b.
\end{array}
\right.
\end{equation}
\begin{theorem}\label{thm:isoG-NM}
For every d.i.  G\"odel algebra ${\bf A}$ and for every d.i.  NM-algebra ${\bf B}$:
\begin{enumerate}
\item The  map $\gamma$ is an isomorphism between ${\bf A}$ and ${\bf G}({\bf NM}^{\pm}(\bf A))$;
\item The map $\eta$ is an isomorphism between ${\bf B}$ and  ${\bf NM^{\pm}(\bf G({\bf B}))}$.
\end{enumerate}
\end{theorem}


\begin{theorem}\label{thm:DINMandG}
A NM-algebra ${\bf B}$ is d.i.
iff there exists a d.i.  G\"odel algebra ${\bf A}$ such that either ${\bf B}\cong{\bf NM}^+({\bf A})$ or ${\bf B}\cong{\bf NM}^-({\bf A})$. 
\end{theorem}

\section{Towards NM-algebras with modal operators}\label{sec:NMOperators}


We first recall the notion of G\"odel algebra with two operators from \cite{FGMR21}. 

\begin{definition}\label{def:GAO}
A {\em G\"odel algebra with operators} ({\em GAO} for short) is a triple $({\bf A}, \Box, \lozenge)$ where ${\bf A}$ is a G\"odel algebra, $\Box$ and $\lozenge$ are unary operators on $A$ satisfying the following equations:
$$
\begin{array}{ll}
(\Box1)\; \Box\top=\top & (\lozenge1)\; \lozenge\bot=\bot\\
(\Box2)\; \Box(x\wedge y)=\Box x\wedge \Box y  \mbox{ }  \mbox{ }  \mbox{ }  & (\lozenge2)\; \lozenge(x\vee y)=\lozenge x\vee \lozenge y.
\end{array}
$$
%
%
\end{definition}
The class of G\"odel algebras with operators is a variety that we denote by $\mathbb{GAO}$. 
{It is easy to check  that every ${\bf A} \in \mathbb{GAO}$ satisfies the following equation (algebraic counterpart of the well-known axiom K in modal logic) and inequations (monotonicity conditions for $\Box$ and $\Diamond$):
\begin{description}
\item [(K)] $\Box (x \to y) \to (\Box x \to \Box y) = \top$.
\item [(Mon)] If $ x \leq y$ then $\Box x \leq \Box y$ and $\Diamond x \leq \Diamond y$.
\end{description}
}

In what follows we will present first steps towards expanding nilpotent minimum algebras with modal operators. In particular we will show that the  construction we recalled in Section \ref{sec:intro} involving d.i.  G\"odel and NM-algebras extends to the case of their modal expansions. 

\subsection{The case of NM$^+$-algebras with operators}\label{secModalNM+}
Let us start by considering a 
GAO $({\bf A}, \Box, \lozenge)$ such that ${\bf A}$ is d.i.  and $\Box$ satisfies the following equation:
\begin{description}
\item[(N1)] $\Box\bot=\bot$.
\end{description}
Let us now define the operators $\minbox$ and $\mindiamond$ on $NM^+({\bf A})$ as follows: for all $(a^-, a^+)\in NM^+({\bf A})$,
\begin{description}
\item[(-)] $\minbox(a^-, a^+)=(\lozenge a^-, \Box a^+)$;
\item[(-)] $\mindiamond(a^-, a^+)=(\Box a^-, \lozenge a^+)$. 
\end{description}
\begin{lemma}
For every $({\bf A}, \Box, \lozenge)$ as above, and for all $(a^-, a^+)\in NM^+({\bf A})$, $\minbox(a^-, a^+), \mindiamond(a^-, a^+)\in NM^+({\bf A})$.
\end{lemma}
\begin{proof2}
Let us prove that $\minbox(a^-, a^+)\in NM^+({\bf A})$, that is to say, $\lozenge a^-\wedge \Box a^+=\bot$ in ${\bf A}$. 
Since ${\bf A}$ is d.i. , $(a^-, a^+)\in NM^+({\bf A})$ iff either $a^-=\bot$ or $a^+=\bot$. If $a^-=\bot$, then $\lozenge a^-\wedge \Box a^+=\lozenge \bot\wedge \Box a^+=\bot \wedge \Box a^+=\bot$. Conversely, if $a^+=\bot$, by $(N1)$, $\Box a^+=\bot$, whence again $\lozenge a^-\wedge \Box a^+=\bot$. This settles the claim.
\end{proof2}
Therefore, the operators $\minbox$ and $\mindiamond$ are well-defined.

\begin{definition}\label{def:nm+mod}
A NM$^+$-algebra with  operators (NMAO$^+$) is a system $({\bf B}, \minbox, \mindiamond)$ where ${\bf B}$ 
is a 
NM$^+$-algebra and $\minbox,\mindiamond: B\to B$ satisfy the following conditions:
$$
\begin{array}{ll}
(\minbox1)\; \minbox\top=\top; &(\mindiamond1) \; \mindiamond\bot=\bot\\
(\minbox2)\; \minbox(x\wedge y)=\minbox x\wedge \minbox y;  \hspace*{0.5cm} & (\mindiamond2)\; \mindiamond(x\vee y)=\mindiamond x\vee \mindiamond y;\\
(F)\; \minbox f= f;& (\minbox\mbox{-}\mindiamond)\; \mindiamond x=\neg \minbox \neg x. 
\end{array}
$$
\end{definition}
The following proposition collects some properties of NMAO$^+$s. 
\begin{proposition}\label{prop:posnegNM+}
In every NMAO$^+$ $({\bf B}, \minbox, \mindiamond)$ the following properties hold:
\begin{enumerate}
\item $\mindiamond f= f$;
\item If $a\geq f$,  $\minbox a\geq f$; if $a\leq f$, $\minbox a\leq f$. The positive and the negative elements of ${\bf B}$ are closed under $\minbox$;
\item If $a\leq f$,  $\mindiamond a\leq f$; if $a\geq f$,  $\mindiamond a\geq f$. The positive and the negative elements of ${\bf B}$ are closed under $\mindiamond$.
\end{enumerate}
\end{proposition}
\begin{proof2}
(1) By $(F)$, $(\minbox\mbox{-}\mindiamond)$ and the fact that $f$ is the negation fixpoint, $\mindiamond f=\neg\minbox\neg f=\neg \minbox f= \neg f=f$.

(2) If $a\geq f$ immediately follows form $(F)$ and the monotonicity property of $\minbox$. If $a\leq f$, then $a\wedge f=a$, whence $\minbox a=\minbox(a\wedge f)=\minbox a\wedge \minbox f$. Thus,  $\minbox a\leq \minbox f=f$.

(3) If $a\leq f$ the claim follows from (1) and the monotonicity property of $\mindiamond$. Conversely, if $a\geq f$, then $a\vee f=a$, whence $\mindiamond a=\mindiamond (a\vee f)=\mindiamond a\vee \mindiamond f$ which shows that $\mindiamond a \geq \mindiamond f=f$. 
\end{proof2}

{ Now we can show that, if we start with a GAO $({\bf A}, \Box, \lozenge)$ with $\bf A$ being d.i. and such that $\Box \bot = \bot$, then the operators $\minbox$ and $\mindiamond$ endow  ${\bf NM}^+({\bf A})$ with the structure of a NMAO.}
\begin{proposition}\label{prop:NM+1}
For every  GAO $({\bf A}, \Box, \lozenge)$ which satisfies (N1) and such that ${\bf A}$ is d.i. , $({\bf NM}^+({\bf A}), \minbox, \mindiamond)$ is a NMAO$^+$. Further, ${\bf NM}^+({\bf A})$ is d.i.  as NM-algebra. 
\end{proposition}
\begin{proof2}
It suffices to prove that $\minbox$ and $\mindiamond$ satisfy the properties of Definition \ref{def:nm+mod}. First of all, recall that the top element of ${\bf NM}^+({\bf A})$ is the pair $(\bot,\top)$. Thus, $\minbox(\bot,\top)=(\lozenge\bot, \Box\top)=(\bot, \top)$ and $(\minbox1)$ holds. Analogously, it is easy to see that $(\mindiamond 1)$ holds as well.

As for $(\minbox2)$, let $(a^-, a^+), (b^-, b^+)\in NM^+({\bf A})$. Then, recalling the definition of $\wedge$ in ${\bf NM}^+({\bf A})$,  $\minbox((a^-, a^+)\wedge (b^-, b^+))=\minbox(a^-\vee b^-, a^+\wedge a^+)=(\lozenge(a^-\vee b^-), \Box(a^+\wedge a^+))=(\lozenge a^-\vee \lozenge b^-, \Box a^+ \wedge \Box b^+)=\minbox(a^-, a^+)\wedge \minbox(b^-, b^+)$. In a similar way one can prove that $(\mindiamond2)$ also holds.  

Now, let $(a^-, a^+)\in NM^+({\bf A})$. Then, $\neg\minbox\neg(a^-, a^+) = \neg\minbox(a^+, a^-) = \neg(\lozenge a^+, \Box a^-) = (\Box a^-, \lozenge a^+) = \mindiamond(a^-, a^+)$ and $(\minbox\mbox{-}\mindiamond)$ holds. 

Let us finally prove $(F)$. The negation fixpoint of ${\bf NM}^+({\bf A})$ is $(\bot, \bot)$ and hence $\minbox(\bot,\bot)=(\lozenge\bot, \Box\bot)$. Since $\lozenge\bot= \bot$ and, by (N1), $\Box\bot= \bot$, by the order relation of ${\bf NM}^+({\bf A})$, $\minbox(\bot,\bot)= (\bot, \bot)$.

That ${\bf NM}^+({\bf A})$ is a d.i.  NM$^+$-algebra follows directly from Theorem \ref{thm:DINMandG} together with the hypothesis that ${\bf A}$ is d.i. 
\end{proof2}

Conversely, let us start with a NMAO$^+$ of the form $({\bf B}, \minbox,\mindiamond)$ where ${\bf B}$ is d.i.  Then, let us define ${\bf G}({\bf B})$ as in Proposition \ref{prop:NMtoGodel} and unary operators $\Box, \lozenge: G({\bf B})\to G({\bf B})$ as: $\Box x=\minbox x\ast \minbox x$ and $\lozenge x=\mindiamond x\ast \mindiamond x$. Then the following holds.
\begin{proposition}\label{prop:NM+2}
For every NMAO$^+$ $({\bf B}, \minbox,\mindiamond)$ where ${\bf B}$ is d.i., $({\bf G}({\bf B}), \Box, \lozenge)$ is a GAO in which ${\bf G}({\bf B})$ is d.i.  and $\Box$ satisfies (N1).  
\end{proposition}
\begin{proof2}
First of all ${\bf G}({\bf B})$ is d.i.  because of Proposition \ref{prop:NMtoGodel}. Moreover (N1) easily holds because, in ${\bf B}$, $b\ast b=\bot$ iff $b\leq f$ and for all such $b$'s, $\minbox b\leq f$ because of Proposition \ref{prop:posnegNM+} (2). Therefore, $\Box \bot=\minbox \bot\ast \minbox \bot=\bot$. 
 Let hence prove that ${\bf G}({\bf B})$ is a GAO. The equations ($\Box1$) and ($\lozenge1$) are easily satisfied. As for the other equations,  let $x,y\in G({\bf B})$, that is, there are $b_x, b_y\in B$ such that $x=b_x\ast b_x$ and $y=b_y\ast b_y$. Let us prove ($\Box2$) distinguishing the following cases:
\vspace{.1cm}

(-) If $b_x, b_y> f$,  $x=b_x\ast b_x=b_x\wedge b_x=b_x$ and $y=b_y\ast b_y=b_y\wedge b_y=b_y$. Thus the claim follows from ($\minbox2$).
\vspace{.1cm}

(-) If $b_x, b_y\leq f$, $x=b_x\ast b_x=\bot=b_y\ast b_y=y$ and the claim trivially holds.
\vspace{.1cm}

(-) If $b_x>f$ and $b_y\leq f$ or $b_x\leq f$ and $b_y> f$, then the claims follows from the previous cases and the observation that, since ${\bf B}$ is d.i. , for all $b_1>f$ and $b_2\leq f$, $b_1\wedge b_2=b_2$.
\vspace{.1cm}

 Finally, ($\lozenge2$) also holds arguing as above. For instance,  if $b_x>f$ and $b_y\leq f$, then $x=b_x$ and $y=\bot$, $\lozenge(x\vee y)=\lozenge x=\lozenge x \vee \lozenge y$ because, by ($\lozenge1$), $\lozenge \bot=\bot$. 
\end{proof2}

\begin{theorem}\label{thm:Gao-NM+}
For every GAO $({\bf A}, \Box,\lozenge)$ which satisfies (N1) and such that ${\bf A}$ is d.i.  and for every NMAO$^+$ $({\bf B}, \minbox, \mindiamond)$ where ${\bf B}$ is d.i., the following claims hold:
\begin{enumerate}
\item  $({\bf A},\Box, \lozenge)\cong ({\bf G}({\bf NM}^+(\bf A)),\Box, \lozenge)$;
\item $({\bf B}, \minbox, \mindiamond)\cong ({\bf NM}^+({\bf G}({\bf B})), \minbox, \mindiamond)$. 
\end{enumerate}
\end{theorem}
\begin{proof2}
In the light of Propositions \ref{prop:NM+1}, \ref{prop:NM+2} and Theorem \ref{thm:isoG-NM}, it suffices to prove that the maps $\gamma$ and $\eta$ defined in (\ref{gamma}) and (\ref{eta}) respectively, preserve the modal operators.  
\vspace{.1cm}

\noindent (1) Let $a\in A$. If $a=\bot$ the claim trivially follows from (N1). Thus, let $a>\bot$ and assume $\Box a>\bot$. Then, $\gamma(\Box a)=(\bot,\Box a)\ast(\bot, \Box a)=(\lozenge\bot,\Box a)\ast(\lozenge\bot, \Box a)=((\Box a\vee \Box a)\to(\lozenge \bot\vee\lozenge \bot), \Box a\wedge\Box a)=(\neg \Box a, \Box a)=(\bot, \Box a)=(\lozenge \bot, \Box a)=(\lozenge \neg a, \Box a)=\minbox(\neg a, a)=\minbox(a\to \bot, a)=\minbox((a\vee a)\to (\bot\vee \bot), a\wedge a)=\minbox((\bot, a)\ast (\bot, a))=\Box(\gamma(a))$, where in the 5th and the 7th equalities we used Lemma \ref{GodProp1} together with the fact that $\Box a>\bot$. If $\Box a=\bot$, $\Box(\gamma(a))=\Box((\bot, a)\ast (\bot, a))=\minbox((a\vee a)\to(\bot\vee \bot), a\wedge a)=\minbox(\neg a, a)=\minbox(\bot, a)=(\lozenge \bot, \Box a)=(\bot, \bot)=\gamma(\bot)=\gamma(\Box a)$.

As for the $\lozenge$ let again $a\in A$, $a>\bot$, $\lozenge a>\bot$. Then  $\gamma(\lozenge(a))=(\bot, \lozenge a)\ast(\bot, \lozenge a)=(\lozenge a\to \bot, \lozenge a)=(\neg \lozenge a, \lozenge a)=(\bot, \lozenge a)=(\Box \bot, \lozenge a)=(\Box \neg a, \lozenge a)=\mindiamond((a\vee a)\to(\bot\vee \bot), a\wedge a)=\lozenge((\bot, a)\ast (\bot, a))=\lozenge(\gamma(a))$. Again we used  Lemma \ref{GodProp1} together with the fact that $\lozenge a>\bot$, and (N1). The case $\lozenge a=\bot$ is analogous to the above and omitted. 
\vspace{.1cm}

\noindent (2) Consider a $b\in B$. If $b=f$, $\eta(\minbox f)=\eta(f)=f=\minbox(f)$. Then, let us take into account the following cases:
\vspace{.1cm}

\noindent (-) $b>f$ and hence $b>\neg b$. Thus, $\eta(\minbox b)=(\bot, \minbox b \ast \minbox b)$. The positive elements of ${\bf B}$ are closed under $\minbox$ (Proposition \ref{prop:posnegNM+} (2)), hence $(\bot, \minbox b \ast \minbox b)=(\bot, \minbox b)=(\mindiamond \bot, \minbox b)=\minbox(\eta(b))$. 
\vspace{.1cm}

\noindent (-) $b<f$ and hence $b<\neg b$. In this case $\eta(\minbox b)=(\neg \minbox b, \bot)=(\mindiamond\neg b, \bot)$. By Proposition \ref{prop:posnegNM+} (3), $\mindiamond\neg b>f$ and hence $\mindiamond\neg b=\mindiamond\neg b\wedge\mindiamond\neg b=\mindiamond\neg b\ast \mindiamond\neg b=\lozenge \neg b$. Thus, $\eta(\minbox b)=(\lozenge \neg b, \bot)=(\lozenge \neg b, \Box \bot)=\minbox\eta(b)$.  

As for the $\mindiamond$, the proof is similar. Let us sketch, for example, the case $b<f$. Again, $\neg b>f$ and $\mindiamond b\leq f$, whence $\eta(\mindiamond b)=(\neg\lozenge b, \bot)=(\minbox \neg b, \bot)=(\Box \neg b, \lozenge \bot)=\mindiamond (\neg b,\bot)=\mindiamond \eta(b)$. 

The claim is hence settled. 
\end{proof2}

\begin{example}\label{ex:Gao-NM+}
Let ${\bf A}$ be the G\"odel algebra whose Hasse diagram is depicted in both (a) and (b) in Figure \ref{fig6} by solid lines,  and let $\Box, \lozenge:A\to A$ be as defined by the dashed arrows in (a) and (b) respectively.  
Notice that $({\bf A}, \Box, \lozenge)$ satisfies (N1).

\begin{figure}
\begin{center}
\begin{tikzpicture}

  \node [label=below:{${\bot}$}, label=below:{} ] (n1)  {} ;
  \node [above  of=n1,label=right:{${a}$}, label=below:{ }] (n2)  {} ;
  \node [above left of=n2,label=left:{${b}$},label=below:{}] (n3) {} ;
  \node [above right of=n2,label=right:{${c}$},label=below:{}] (n4) {} ;
  \node [above left of=n4,label=above:{${\top}$},label=below:{}] (n6) {} ;

  \draw  (n1) -- (n2);
  \draw (n2) -- (n3);
  \draw  (n3) -- (n6);
  \draw  (n2) -- (n4);
  \draw  (n4) -- (n6);
  \draw [fill] (n1) circle [radius=.5mm];
  \draw [fill] (n2) circle [radius=.5mm];
  \draw [fill] (n3) circle [radius=.5mm];
  \draw [fill] (n4) circle [radius=.5mm];
  \draw [fill] (n6) circle [radius=.5mm];

  \draw [->,dashed] (n2) arc [radius=3mm, start angle=90, end angle= 400](n2);
  \draw [->,dashed] (n4) edge[bend right] (n3);
   \draw [->,dashed] (n3) edge[bend right] (n4);
  \draw [->,dashed] (n1) arc [radius=3mm, start angle=90, end angle= 400]  (n1);
  \draw [->,dashed] (n6) arc [radius=3mm, start angle=270, end angle=610 ]  (n6);
\end{tikzpicture}
\begin{tikzpicture}

  \node [label=below:{${\bot}$}, label=below:{} ] (n1)  {} ;
  \node [above  of=n1,label=right:{${a}$}, label=below:{ }] (n2)  {} ;
  \node [above left of=n2,label=left:{${b}$},label=below:{}] (n3) {} ;
  \node [above right of=n2,label=right:{${c}$},label=below:{}] (n4) {} ;
  \node [above left of=n4,label=above:{${\top}$},label=below:{}] (n6) {} ;

  \draw  (n1) -- (n2);
  \draw (n2) -- (n3);
  \draw  (n3) -- (n6);
  \draw  (n2) -- (n4);
  \draw  (n4) -- (n6);
  \draw [fill] (n1) circle [radius=.5mm];
  \draw [fill] (n2) circle [radius=.5mm];
  \draw [fill] (n3) circle [radius=.5mm];
  \draw [fill] (n4) circle [radius=.5mm];
  \draw [fill] (n6) circle [radius=.5mm];

  \draw [->,dashed] (n2) edge[bend left] (n1);
  \draw [->,dashed] (n3) edge[bend right] (n2);
  \draw [->,dashed] (n4) edge[bend left] (n2);
  \draw [->,dashed] (n1) arc [radius=3mm, start angle=90, end angle= 400]  (n1);
  \draw [->,dashed] (n6) edge  (n2);
\end{tikzpicture}
\begin{tikzpicture}

  \node [label=below:{${\bot}$}, label=below:{} ] (n1)  {} ;
  \node [above left of=n1,label=left:{${x}$}, label=below:{ }] (n2)  {} ;
  \node [above right of=n1,label=right:{${y}$},label=below:{}] (n3) {} ;
  \node [above left of=n3,label=right:{${z}$},label=below:{}] (n4) {} ;
  \node [above  of=n4,label=left:{${f}$},label=below:{}] (n6) {} ;
   \node [above  of=n6,label=left:{${k}$},label=below:{}] (n7) {} ;
      \node [above left of=n7,label=left:{${t}$},label=below:{}] (n8) {} ;
         \node [above right of=n7,label=right:{${s}$},label=below:{}] (n9) {} ;
            \node [above right of=n8,label=above:{${\top}$},label=below:{}] (n10) {} ;

  \draw  (n1) -- (n2);
  \draw (n1) -- (n3);
  \draw  (n3) -- (n4);
  \draw  (n2) -- (n4);
  \draw  (n4) -- (n6);
  \draw  (n6) -- (n7);
  \draw (n7)--(n8);
  \draw (n7)--(n9);
  \draw (n8)--(n10);
  \draw (n9)--(n10);
  \draw [fill] (n1) circle [radius=.5mm];
  \draw [fill] (n2) circle [radius=.5mm];
  \draw [fill] (n3) circle [radius=.5mm];
  \draw [fill] (n4) circle [radius=.5mm];
  \draw [fill] (n6) circle [radius=.5mm];
   \draw [fill] (n7) circle [radius=.5mm];
    \draw [fill] (n8) circle [radius=.5mm];
     \draw [fill] (n9) circle [radius=.5mm];
      \draw [fill] (n10) circle [radius=.5mm];
  \draw [->,dashed] (n1) edge (n4);
    \draw [->,dashed] (n2) edge [bend left]  (n4);
        \draw [->,dashed] (n3) edge [bend right]  (n4);
          \draw [->,dashed] (n4) edge [bend right]  (n6);
  \draw [->,dashed] (n10) arc [radius=3mm, start angle=270, end angle= 600]  (n10);
  \draw [->,dashed] (n6) arc [radius=3mm, start angle=0, end angle= 320]  (n6);
   \draw [->,dashed] (n7) arc [radius=3mm, start angle=0, end angle= 320]  (n7);
   \draw [->,dashed] (n8) edge [bend right]  (n9);
   \draw [->,dashed] (n9) edge [bend right]  (n8);
\end{tikzpicture}
\begin{tikzpicture}

  \node [label=below:{${\bot}$}, label=below:{} ] (n1)  {} ;
  \node [above left of=n1,label=left:{${x}$}, label=below:{ }] (n2)  {} ;
  \node [above right of=n1,label=right:{${y}$},label=below:{}] (n3) {} ;
  \node [above left of=n3,label=right:{${z}$},label=below:{}] (n4) {} ;
  \node [above  of=n4,label=left:{${f}$},label=below:{}] (n6) {} ;
   \node [above  of=n6,label=left:{${k}$},label=below:{}] (n7) {} ;
      \node [above left of=n7,label=left:{${t}$},label=below:{}] (n8) {} ;
         \node [above right of=n7,label=right:{${s}$},label=below:{}] (n9) {} ;
            \node [above right of=n8,label=above:{${\top}$},label=below:{}] (n10) {} ;

  \draw  (n1) -- (n2);
  \draw (n1) -- (n3);
  \draw  (n3) -- (n4);
  \draw  (n2) -- (n4);
  \draw  (n4) -- (n6);
  \draw  (n6) -- (n7);
  \draw (n7)--(n8);
  \draw (n7)--(n9);
  \draw (n8)--(n10);
  \draw (n9)--(n10);
  \draw [fill] (n1) circle [radius=.5mm];
  \draw [fill] (n2) circle [radius=.5mm];
  \draw [fill] (n3) circle [radius=.5mm];
  \draw [fill] (n4) circle [radius=.5mm];
  \draw [fill] (n6) circle [radius=.5mm];
   \draw [fill] (n7) circle [radius=.5mm];
    \draw [fill] (n8) circle [radius=.5mm];
     \draw [fill] (n9) circle [radius=.5mm];
      \draw [fill] (n10) circle [radius=.5mm];
  \draw [->,dashed] (n1) arc [radius=3mm, start angle=90, end angle= 400] (n1);
    \draw [->,dashed] (n2) edge [bend right]  (n3);
        \draw [->,dashed] (n3) edge [bend right]  (n2);
          \draw [->,dashed] (n4) arc [radius=3mm, start angle=0, end angle= 320](n4);
  \draw [->,dashed] (n10) edge (n7);
  \draw [->,dashed] (n6) arc [radius=3mm, start angle=0, end angle= 320]  (n6);
   \draw [->,dashed] (n7) edge [bend left] (n6);
   \draw [->,dashed] (n8) edge [bend right]  (n7);
   \draw [->,dashed] (n9) edge [bend left]  (n7);
%
\end{tikzpicture} \vspace{.2cm}

 \hfill (a) \hfill (b) \hfill (c)  \hfill (d) \hfill \mbox{}
\end{center}
\caption{{\small \small A d.i. G\"odel algebra with a $\Box$ operator satisfying (N1) (a) and a $\lozenge$ operator (b), and the corresponding d.i. NM$^+$-algebra with the operators $\minbox$ (c) and $\minbox$ (d).}}
\label{fig6}
\end{figure}

The Hasse diagram of the algebra ${\bf NM}^+({\bf A})$ is the one depicted in (c) and (d) in Figure \ref{fig6} with solid lines,  where: 
$\bot=(\top, \bot)$, $x=(c,\bot)$, $y=(b, \bot)$, $z=(a,\bot)$, $f=(\bot,\bot)$, $k=(\bot, a)$, $t=(\bot, b)$, $s=(\bot,c)$, $\top=(\bot,\top)$. The operators $\minbox$ and $\mindiamond$ on $NM^+({\bf A})$ defined as above
correspond to the dashed arrows depicted in (c) and (d) of Figure \ref{fig6} respectively. 
For instance, $\minbox z=\minbox(a,\bot)=(\lozenge a, \Box \bot)=(\bot,\bot)=f$ and $\mindiamond t=\mindiamond (\bot, b)=(\Box\bot,\lozenge b)=(\bot, a)=k$. 
\end{example}
\subsection{The case of NM$^-$-algebras with operators}\label{secModalNM-}
Let $({\bf A}, \Box, \lozenge)$ be a GAO in which ${\bf A}$ is d.i.  and satisfying the condition (N1) and the following additional ones:
$$
\begin{array}{ll}
(SM\Box)\, \mbox{ If }a>\bot\mbox{, then }\Box a>\bot; & \hspace*{0.3cm}
(SM\lozenge)\, \mbox{ If }a>\bot\mbox{, then }\lozenge a>\bot.
\end{array}
$$
Let us define $\minbox$ and $\mindiamond$ on $NM^-({\bf A})$ as above: for all $(a^-, a^+)\in NM^-({\bf A})$, $\minbox(a^-, a^+)=(\lozenge a^-, \Box a^+)$ and $\mindiamond(a^-, a^+)=(\Box a^-, \lozenge a^+)$.
\begin{lemma}
For every $({\bf A}, \Box, \lozenge)$ as above, and for all $(a^-, a^+)\in NM^-({\bf A})$, $\minbox(a^-, a^+), \mindiamond(a^-, a^+)\in NM^-({\bf A})$.
\end{lemma}
\begin{proof2}
Let us prove the claim for the case of $\minbox$ and in particular that for all $(a^-, a^+)$, if it satisfies $(a^-\wedge a^+)\vee \neg (a^-\vee a^+)=\bot$, then $(\lozenge a^-\wedge \Box a^+)\vee \neg (\lozenge a^-\vee \Box a^+)=\bot$. Since ${\bf A}$ is directly indecomposable, as G\"odel algebra, then so is $NM^-({\bf A})$ as NM-algebra. Therefore, $(a^-\wedge a^+)\vee \neg (a^-\vee a^+)=\bot$ iff either $a^-=\bot$ and $a^+>\bot$, or $a^->\bot$ and $a^+=\bot$.  
If the former is the case, then $\lozenge a^-=\bot$ by $(\lozenge1)$ and $\Box a^+>\bot$  by (SM$\Box$). 
If the latter is the case, similarly, $\lozenge a^->\bot$ by (SM$\lozenge$) and $\Box a^+=\bot$ because of (N1). 
Thus in both cases  $(\lozenge a^-\wedge \Box a^+)\vee \neg (\lozenge a^-\vee \Box a^+)=\bot$, which settles the claim.
\end{proof2}

For the following, recall the equations of Definition \ref{def:nm+mod}.
\begin{definition}\label{def:nm+mod}
An NM$^-$-algebra with  operators (NMAO$^-$) is a system $({\bf B}, \Box, \lozenge)$ where ${\bf B}=(B, \ast, \wedge, \neg, \bot)$ be a 
NM$^-$-algebra, and $\Box,\lozenge: B\to B$ satisfy the following conditions: $(\minbox1)$; $(\minbox2)$; $(\minbox\mbox{-}\mindiamond)$; $(\mindiamond1)$; $(\mindiamond2)$; and

$
(P):\; \minbox(x\vee \neg x)=\minbox x\vee \neg \minbox x$; \hspace*{0.2cm}  
$(N):\; \mbox{ if }x\leq \neg x\mbox{, then }\mindiamond x \leq  \neg \mindiamond x. $

%
\end{definition}

\begin{proposition}\label{prop:NM-1}
The following properties hold in every NMAO$^-$:
\begin{enumerate}
\item The positive and the negative elements of ${\bf B}$ are closed under $\minbox$;
\item The positive and the negative elements of ${\bf B}$ are closed under $\mindiamond$.
\end{enumerate}
\end{proposition}

\begin{proof2}
(1) Suppose $x\geq \neg x$. Then $x\vee\neg x=x$ and hence, by (P), $\minbox x=\minbox(x\vee\neg x)=\minbox x\vee \neg\minbox x$. Thus, $\minbox x\geq \neg \minbox x$. Conversely, if $x\leq \neg x$, by (P), $\minbox \neg x=\minbox(x\vee \neg x)=\minbox x\vee \neg \minbox x$, whence $\minbox x\leq \minbox\neg x$.

(2) The first part of the claim follows from (N), the second by the order-reversing property of $\neg$. 
\end{proof2}

 Now we prove that, if {\bf A} is a d.i. G\"odel algebra and $({\bf A}, \Box, \lozenge)$ is a GAO satisfying (N1), (SM$\Box$) and (SM$\lozenge$), then  ${\bf NM}^-({\bf A})$ is a NM$^-$-algebra with operators.
\begin{proposition}\label{prop:NM-1}
For any GAO $({\bf A}, \Box, \lozenge)$ satisfying (N1), (SM$\Box$) and (SM$\lozenge$) such that ${\bf A}$ is d.i., $({\bf NM}^-({\bf A}), \minbox, \mindiamond)$ is a NMAO$^-$ and ${\bf NM}^-({\bf A})$ is d.i.  as NM-algebra. 
\end{proposition}
\begin{proof2}
The equations ($\minbox1$), ($\minbox2$), ($\mindiamond1$), ($\mindiamond2$) and ($\minbox\mbox{-}\mindiamond$) hold with the same proof  of Proposition \ref{prop:NM+1}. Let us hence prove (P) and (N). 

If $(a^-, a^+)\in NM^-({\bf A})$, then either $a^-=\bot$ or $a^+=\bot$. We assume that $a^-=\bot$ without loss of generality (the case $a^+=\bot$ is symmetric and omitted). Then, $\minbox((a^-, a^+)\vee\neg(a^-, a^+))=\minbox(a^-\wedge a^+, a^+\vee a^-)=\minbox(\bot, a^+)=(\lozenge\bot, \Box a^+)=(\bot, \Box a^+)$. On the other hand, $\minbox(a^-, a^+)\vee \neg \minbox(a^-, a^+)=(\lozenge \bot, \Box a^+)\vee \neg (\lozenge \bot, \Box a^+)=(\bot, \Box a^+)\vee (\Box a^+, \bot)=(\bot\wedge \Box a^+, \bot\vee \Box a^+)=(\bot, \Box a^+)$.  Thus (P) holds.

As for (N), if $(a^-, a^+)\leq \neg (a^-, a^+)$, then $a^-\geq \bot$ and $a^+=\bot$. Therefore,  $\mindiamond(a^-, a^+)=(\Box a^-, \lozenge a^+)=(\Box a^-, \bot)$. On the other hand $\neg \boxminus(a^-, a^+)=\mindiamond(a^+, a^-)=(\Box a^+, \lozenge a^-)=(\bot, \lozenge a^-)$ and $(\Box a^-, \bot)\leq (\bot, \lozenge a^-)$. 
\end{proof2}

And this is the converse direction, from  NMAO$^-$s to GAOs.
\begin{proposition}\label{prop:NM-2}
For every NMAO$^-$ $({\bf B}, \minbox,\mindiamond)$ where ${\bf B}$ is d.i. , $({\bf G}({\bf B}), \Box, \lozenge)$ is a GAO where ${\bf G}({\bf B})$ is d.i.  and $\Box$ satisfies (N1), (SM$\Box$) and (SM$\lozenge$).  
\end{proposition}
\begin{proof2}
For every d.i.  NM$^-$-algebra ${\bf B}$, Proposition \ref{prop:NMtoGodel} tells us that ${\bf G}({\bf B})$ is a d.i.  G\"odel algebra.  The  proofs of ($\Box1$), ($\Box2$), ($\lozenge 1$), ($\lozenge 2$) and (N1) are as in Proposition \ref{prop:NM+2}, hence let us prove (SM$\Box$) and (SM$\lozenge$). 

Take  $a>\bot$ in ${\bf G}({\bf B})$ and  let $b\in B$ such that $a=b\ast b>\bot$. Thus $b>\neg b$ and  $a=b\wedge b=b$. It follows that $\Box a=\minbox a\ast \minbox a=\minbox b\ast \minbox b$. By Prop. \ref{prop:NM-1}, $\minbox b>\neg \minbox b$, whence $\Box a=\minbox b\ast \minbox b=\minbox b>\neg \minbox b$ and hence, in ${\bf G}({\bf B})$, $\Box a>\bot$. Similarly, $\lozenge a=\lozenge b\ast\lozenge b=\lozenge b$ (again by Prop. \ref{prop:NM-1}) whence, in ${\bf G}({\bf B})$, $\lozenge a>\bot$. 
\end{proof2}
Essentially the same proof of Theorem \ref{thm:Gao-NM+}, 
together with the above Propositions \ref{prop:NM-1} and \ref{prop:NM-2}, proves the following representation theorem.
\begin{theorem}\label{thm:Gao-NM-}
For all GAO $({\bf A}, \Box,\lozenge)$ which satisfies (N1), (SM$\Box$) and (SM$\lozenge$) and  ${\bf A}$ is d.i.  and for all NMAO$^-$ $({\bf B}, \minbox, \mindiamond)$ where ${\bf B}$ is d.i.,  we have that $({\bf A},\Box, \lozenge)\cong ({\bf G}({\bf NM}^-(\bf A)),\Box, \lozenge)$ and 
 $({\bf B}, \minbox, \mindiamond)\cong ({\bf NM}^-({\bf G}({\bf B})), \minbox, \mindiamond)$. 
\end{theorem}
Figure \ref{fig7}  shows the effect of the construction NM$^-$ to the same GAO we already discussed in Example \ref{ex:Gao-NM+} and whose G\"odel algebra reduct is directly indecomposable. 

\begin{figure}
\begin{center}
\begin{tikzpicture}

  \node [label=below:{${\bot}$}, label=below:{} ] (n1)  {} ;
  \node [above  of=n1,label=right:{${a}$}, label=below:{ }] (n2)  {} ;
  \node [above left of=n2,label=left:{${b}$},label=below:{}] (n3) {} ;
  \node [above right of=n2,label=right:{${c}$},label=below:{}] (n4) {} ;
  \node [above left of=n4,label=above:{${\top}$},label=below:{}] (n6) {} ;

  \draw  (n1) -- (n2);
  \draw (n2) -- (n3);
  \draw  (n3) -- (n6);
  \draw  (n2) -- (n4);
  \draw  (n4) -- (n6);
  \draw [fill] (n1) circle [radius=.5mm];
  \draw [fill] (n2) circle [radius=.5mm];
  \draw [fill] (n3) circle [radius=.5mm];
  \draw [fill] (n4) circle [radius=.5mm];
  \draw [fill] (n6) circle [radius=.5mm];

  \draw [->,dashed] (n2) arc [radius=3mm, start angle=90, end angle= 400](n2);
  \draw [->,dashed] (n4) edge[bend right] (n3);
   \draw [->,dashed] (n3) edge[bend right] (n4);
  \draw [->,dashed] (n1) arc [radius=3mm, start angle=90, end angle= 400]  (n1);
  \draw [->,dashed] (n6) arc [radius=3mm, start angle=270, end angle=610 ]  (n6);
\end{tikzpicture}
\begin{tikzpicture}

  \node [label=below:{${\bot}$}, label=below:{} ] (n1)  {} ;
  \node [above  of=n1,label=right:{${a}$}, label=below:{ }] (n2)  {} ;
  \node [above left of=n2,label=left:{${b}$},label=below:{}] (n3) {} ;
  \node [above right of=n2,label=right:{${c}$},label=below:{}] (n4) {} ;
  \node [above left of=n4,label=above:{${\top}$},label=below:{}] (n6) {} ;

  \draw  (n1) -- (n2);
  \draw (n2) -- (n3);
  \draw  (n3) -- (n6);
  \draw  (n2) -- (n4);
  \draw  (n4) -- (n6);
  \draw [fill] (n1) circle [radius=.5mm];
  \draw [fill] (n2) circle [radius=.5mm];
  \draw [fill] (n3) circle [radius=.5mm];
  \draw [fill] (n4) circle [radius=.5mm];
  \draw [fill] (n6) circle [radius=.5mm];

  \draw [->,dashed] (n2) arc [radius=3mm, start angle=90, end angle= 400]  (n2);
  \draw [->,dashed] (n3) edge[bend right] (n2);
  \draw [->,dashed] (n4) edge[bend left] (n2);
  \draw [->,dashed] (n1) arc [radius=3mm, start angle=90, end angle= 400]  (n1);
  \draw [->,dashed] (n6) edge  (n2);
\end{tikzpicture}
\begin{tikzpicture}

  \node [label=below:{${\bot}$}, label=below:{} ] (n1)  {} ;
  \node [above left of=n1,label=left:{${x}$}, label=below:{ }] (n2)  {} ;
  \node [above right of=n1,label=right:{${y}$},label=below:{}] (n3) {} ;
  \node [above left of=n3,label=right:{${z}$},label=below:{}] (n4) {} ;
   \node [above  of=n4,label=left:{${k}$},label=below:{}] (n7) {} ;
      \node [above left of=n7,label=left:{${t}$},label=below:{}] (n8) {} ;
         \node [above right of=n7,label=right:{${s}$},label=below:{}] (n9) {} ;
            \node [above right of=n8,label=above:{${\top}$},label=below:{}] (n10) {} ;

  \draw  (n1) -- (n2);
  \draw (n1) -- (n3);
  \draw  (n3) -- (n4);
  \draw  (n2) -- (n4);
  \draw  (n4) --  (n7);
  \draw (n7)--(n8);
  \draw (n7)--(n9);
  \draw (n8)--(n10);
  \draw (n9)--(n10);
  \draw [fill] (n1) circle [radius=.5mm];
  \draw [fill] (n2) circle [radius=.5mm];
  \draw [fill] (n3) circle [radius=.5mm];
  \draw [fill] (n4) circle [radius=.5mm];
   \draw [fill] (n7) circle [radius=.5mm];
    \draw [fill] (n8) circle [radius=.5mm];
     \draw [fill] (n9) circle [radius=.5mm];
      \draw [fill] (n10) circle [radius=.5mm];
  \draw [->,dashed] (n1) edge (n4);
    \draw [->,dashed] (n2) edge [bend left]  (n4);
        \draw [->,dashed] (n3) edge [bend right]  (n4);
          \draw [->,dashed] (n4)arc [radius=2.4mm, start angle=270, end angle=610 ]   (n4);
  \draw [->,dashed] (n10) arc [radius=3mm, start angle=270, end angle= 600]  (n10);
   \draw [->,dashed] (n7) arc [radius=2.4mm, start angle=90, end angle= 400] (n7);
   \draw [->,dashed] (n8) edge [bend right]  (n9);
   \draw [->,dashed] (n9) edge [bend right]  (n8);
\end{tikzpicture}
\begin{tikzpicture}

   \node [label=below:{${\bot}$}, label=below:{} ] (n1)  {} ;
  \node [above left of=n1,label=left:{${x}$}, label=below:{ }] (n2)  {} ;
  \node [above right of=n1,label=right:{${y}$},label=below:{}] (n3) {} ;
  \node [above left of=n3,label=right:{${z}$},label=below:{}] (n4) {} ;
   \node [above  of=n4,label=left:{${k}$},label=below:{}] (n7) {} ;
      \node [above left of=n7,label=left:{${t}$},label=below:{}] (n8) {} ;
         \node [above right of=n7,label=right:{${s}$},label=below:{}] (n9) {} ;
            \node [above right of=n8,label=above:{${\top}$},label=below:{}] (n10) {} ;

  \draw  (n1) -- (n2);
  \draw (n1) -- (n3);
  \draw  (n3) -- (n4);
  \draw  (n2) -- (n4);
  \draw  (n4) --  (n7);
  \draw (n7)--(n8);
  \draw (n7)--(n9);
  \draw (n8)--(n10);
  \draw (n9)--(n10);
  \draw [fill] (n1) circle [radius=.5mm];
  \draw [fill] (n2) circle [radius=.5mm];
  \draw [fill] (n3) circle [radius=.5mm];
  \draw [fill] (n4) circle [radius=.5mm];
   \draw [fill] (n7) circle [radius=.5mm];
    \draw [fill] (n8) circle [radius=.5mm];
     \draw [fill] (n9) circle [radius=.5mm];
      \draw [fill] (n10) circle [radius=.5mm];
  \draw [->,dashed] (n1) arc [radius=3mm, start angle=90, end angle= 400] (n1);
    \draw [->,dashed] (n2) edge [bend right]  (n3);
        \draw [->,dashed] (n3) edge [bend right]  (n2);
          \draw [->,dashed] (n4) arc [radius=2.4mm, start angle=270, end angle=610 ]   (n4);
  \draw [->,dashed] (n10) edge (n7);
   \draw [->,dashed] (n7) arc [radius=2.4mm, start angle=90, end angle= 400]    (n7);
   \draw [->,dashed] (n8) edge [bend right]  (n7);
   \draw [->,dashed] (n9) edge [bend left]  (n7);
%
\end{tikzpicture}  \vspace{.2cm}

 \hfill (a) \hfill (b) \hfill (c)  \hfill (d) \hfill \mbox{}
\end{center}
\caption{{\small \small A d.i. G\"odel algebra endowed with a $\Box$ operator (a)  and a $\lozenge$ operator (b) satisfying (N1), (SM$\Box$) and (SM$\lozenge$), and the corresponding d.i. NM$^-$-algebra with the operators $\minbox$ (c) and $\mindiamond$ (d).}}
\label{fig7}
\end{figure}

\section{Conclusion and future work}\label{secFinal}

In this  paper we have studied the effect of taking connected and disconnected rotations of G\"odel algebras with operators whose G\"odel reduct is directly indecomposable. By doing so, we have introduced the varieties of NM$^+$ and NM$^-$-algebras with modal operators and we have characterized the members of these varieties whose NM-reduct is directly indecomposable.

Although the aim of the present paper is to put forward an algebraic analysis of the modal structures we took into account, in \cite{FGR19,FGMR21} we also gave a description of the dual relational frames that  arise from finite G\"odel algebras with operators ({\em forest frames}). Those are 
the prelinear version of the  models, based on posets, that are a semantics for intuitionistic modal logic (see e.g., \cite{OrlRew,Ale}).  
In our future work, we plan to extend the present approach essentially in thee directions:

1. From the point of view of forest frames, taking into account that  for every finite d.i. G\"odel algebra ${\bf A}$,  ${\bf NM}^+({\bf A})$ and ${\bf NM}^-({\bf A})$ have the same prime filters, we plan to investigate how NMAO$^+$ and NMAO$^-$-algebras relate to forest frames and to extend to the NM-case the isomorphic representation theorem from \cite{FGR19}. 

2. From an algebraic perspective, besides connected and disconnected rotations, NM-algebras can be also seen as twist-structures obtained from G\"odel algebras. It is hence interesting to deepen this latter construction for modal G\"odel algebras and compare the modal NM-algebras obtained in such a way with the modal NM$^+$- and NM$^-$-algebras defined in the present paper, also in light of the general results proved in \cite{RivOno}. 

3. From a more general point of view, NM-algebras can be also seen as a subvariety of Nelson lattices. In \cite{MR2021}, the authors introduce a definition of algebra with operators more general than the considered one in the current paper in the sense that they are not required to satisfy {\bf Mon} (monotony rules for modal operators). It will be interesting to compare this approach to ours.    
\vspace{.1cm}

\noindent{\bf Acknowledgments}. The authors  thank 
  the anonymous referees for their comments. 
Authors acknowledge partial support by the MOSAIC project (EU H2020-MSCA-RISE-2020 Project 101007627). Flaminio and Godo also acknowledge partial support by the Spanish project PID2019-111544GB-C21 funded by  \\ MCIN/AEI/10.13039/501100011033. Menchon acknowledge partial support by argentinean projects PIP 112-20200101301CO (CONICET) and PICT-2019-2019-00882 (ANPCyT). The fourth author wants to acknowledge partial support by the  following argentinean projects: PIP 112-20150100412CO (CONICET) and UBA-CyT-20020190100021BA.

\end{document}